\documentclass{article}
\usepackage{amssymb,amsmath}

\newtheorem{theorem}{Theorem}
\newtheorem{lemma}{Lemma}
\newtheorem{problem}{Problem}
\newcommand{\bt}{\begin{theorem}}
\newcommand{\et}{\end{theorem}}
\newcommand{\bl}{\begin{lemma}}
\newcommand{\el}{\end{lemma}}
\newcommand{\bp}{\begin{problem}}
\newcommand{\ep}{\end{problem}}
\newcommand{\pf}{{\bf Proof}.\ }
\newcommand{\eop}{$\square$\vspace{.8cm}}
\newcommand{\bq}{\begin{eqnarray*}}
\newcommand{\eq}{\end{eqnarray*}}
\newcommand{\be}{\begin{eqnarray}}
\newcommand{\ee}{\end{eqnarray}}
\newcommand{\beq}{\begin{equation}}
\newcommand{\eeq}{\end{equation}}
\newcommand{\benum}{\begin{enumerate}}
\newcommand{\eenum}{\end{enumerate}}
\newcommand{\ba}{\begin{array}}
\newcommand{\ea}{\end{array}}
\newcommand{\Z}{\ensuremath{\mathbf Z}}
\newcommand{\N}{\ensuremath{ \mathbf N }}
\newcommand{\FF}{\ensuremath{\mathcal F}}
\newcommand{\GG}{\ensuremath{\mathcal G}}
\newcommand{\HH}{\ensuremath{\mathcal H}}
\newcommand{\supp}{\mbox{supp}}

\newcommand{\addq}{\ensuremath{\oplus_q}}
\newcommand{\pol}{$\mathcal{F} = \{f_n(q)\}_{n=1}^{\infty}$}
\newcommand{\polg}{$\mathcal{G} = \{g_n(q)\}_{n=1}^{\infty}$}
\newcommand{\tfe}{ the functional equation~(\ref{q:fe})}
\newcommand{\fe}{~(\ref{q:fe})}

\date{}

\begin{document}
\title{A functional equation arising from \\
multiplication of quantum integers}
\author{Melvyn B. Nathanson\thanks{This work was supported
in part by grants from the NSA Mathematical Sciences Program
and the PSC-CUNY Research Award Program.}\\
Department of Mathematics\\
Lehman College (CUNY)\\
Bronx, New York 10468\\
Email: nathansn@alpha.lehman.cuny.edu}
\maketitle

\begin{abstract}
For the quantum integer $[n]_q = 1+q+q^2+\cdots + q^{n-1}$
there is a natural polynomial multiplication
such that $[m]_q\otimes_q [n]_q = [mn]_q$.  This multiplication leads 
to the functional equation $f_m(q)f_n(q^m) = f_{mn}(q)$, defined on a given sequence 
\pol\ of polynomials.
This paper contains various results concerning the construction
and classification of polynomial sequences that satisfy the functional equation,
as well open problems that arise from the functional equation.
\end{abstract}

\section{A polynomial functional equation}
Let \N$ = \{1,2,3,\ldots\}$ denote the set of natural numbers,
and $\N_0 = \N \cup \{0\}$ the set of nonnegative integers.
For $n \in \N$, the polynomial
\[
[n]_q = 1 + q + q^2 + \cdots + q^{n-1}
\]
is called the {\em quantum integer} $n.$
With the usual multiplication of polynomials, however, we observe that
$[m]_q [n]_q \neq [mn]_q$ for all $m\neq 1$ and $n\neq 1$.
We would like to define a polynomial multiplication such that
the product of the quantum integers $[m]_q$ and $[n]_q$
is $[mn]_q$.

Consider polynomials with coefficients in a field.
Let $\mathcal{F} = \{f_n(q)\}_{n=1}^{\infty}$ be a sequence of polynomials.
We define a multiplication operation on the polynomials in $\mathcal{F}$ by
\[
f_m(q)\otimes_q f_n(q) = f_m(q)f_n(q^m).
\]
We want to determine all sequences $\mathcal{F}$
that satisfy the {\em functional equation}
\begin{equation}       \label{q:fe}
f_{mn}(q) = f_m(q)\otimes_q f_n(q) = f_m(q)f_n(q^m)
\eeq
for all $m,n \in \N$.
If the sequence $\mathcal{F}$
is a solution of~(\ref{q:fe}), then the operation $\otimes_q$
is commutative on $\mathcal{F}$ since
\[
f_m(q)\otimes_q f_n(q) = f_{mn}(q) = f_{nm}(q) = f_n(q)\otimes_q f_m(q).
\]
Equivalently,
\beq                 \label{q:mnfe}
f_m(q)f_n(q^m) = f_n(q)f_m(q^n)
\eeq
for all natural numbers $m$ and $n$.\footnote{Note that\fe\
implies~(\ref{q:mnfe}), but not conversely, since the sequence \pol\
with $f_n(q) = 2$ for all $n \in \N$ satisfies~(\ref{q:mnfe})
but not\fe.}

Here are three examples of solutions of \tfe.
First, the constant sequence defined by $f_n(q) = 1$ for 
all $n \in \N$ satisfies~(\ref{q:fe}).

Second, let
\[
f_n(q) = q^{n-1}
\]
for all $n \in \N.$  Then
\[
f_{mn}(q) = q^{mn-1} = q^{m-1}q^{m(n-1)} = f_m(q)f_n\left( q^m\right),
\]
and so the polynomial sequence $\{q^{n-1}\}_{n=1}^{\infty}$
also satisfies~(\ref{q:fe}).

Third, let $f_n(q) = [n]_q$ for all $n \in \N$.  Then
\bq
[m]_q\otimes_q [n]_q
& = & f_m(q)\otimes_q f_n(q)  \\
& = & f_m(q)f_n(q^m)  \\
& = & (1+q+q^2 + \cdots + q^{m-1})(1+q^m+q^{2m}+\cdots + q^{m(n-1)})\\
& = & 1+q + \cdots + q^{m-1} + q^m +q^{m+1}\cdots + q^{mn-1}\\
& = & [mn]_q,
\eq
and so the polynomial sequence $\{[n]_q\}_{n=1}^{\infty}$
of quantum integers satisfies the functional equation~\eqref{q:fe}.

The identity 
\[
[m]_q\otimes_q [n]_q = [mn]_q
\]
is the $q$-series expression of the following additive number theoretic 
identity for sumsets:
\[
\{0,1,2,\ldots,mn-1\} 
= \{0,1,\ldots,m-1\} + \{ 0,m,2m,\ldots, (n-1)m\} \\
\]

This paper investigates the following problem:
\bp          \label{q:p1}
Determine all polynomial 
sequences \pol\ that satisfy the functional equation~(\ref{q:fe}).
\ep

\section{Prime semigroups}
A {\em multiplicative subsemigroup} of the natural numbers, or, simply,
a {\em semigroup}, is a set $S \subseteq \N$
such that $1 \in S$ and if $m \in S$ and $n \in S$, then $mn \in S.$
For example, for any positive integer $n_0,$ the set
$\{1\} \cup \{n \geq n_0\}$ is a semigroup.
If $P$ is a set of prime numbers, then the set $S(P)$
consisting of the positive integers all of whose prime factors belong to $P$
is a multiplicative subsemigroup of \N.  
If $P = \emptyset,$ then $S(P) = \{1\}.$
If $P = \{p\}$ contains only one prime, then
$S(P) = \{p^k:k\in \N_0\}.$
A semigroup of the form $S(P)$, where $P$ is a set of primes, 
will be called a {\em prime semigroup.}

Let $\mathcal{F} = \{f_n(q)\}_{n=1}^{\infty}$
be a sequence of polynomials that satisfies
the functional equation~(\ref{q:fe}).
Since
\[
f_1(q) = f_1(q)\otimes_q f_1(q) = f_1(q)f_1(q),
\]
it follows that $f_1(q) = 1$ or $f_1(q) = 0$.
If $f_1(q) = 0$, then
\[
f_n(q) = f_1(q)\otimes_q f_n(q) = f_1(q)f_n(q) = 0
\]
for all $n \in \N$, and $\mathcal{F}$ is the sequence of zero polynomials.
If $f_n(q) \neq 0$ for some $n$, then $f_1(q) = 1$.

Let \pol\ be any sequence of functions.
The {\em support} of \FF\ is the set
\[
\supp(\FF) = \{n \in \N : f_n(q) \neq 0\}.
\]
The sequence \FF\ is called {\em nonzero} if $f_n(q)\neq 0$ for some $n \in\N$,
or, equivalently, if $\supp(\FF)\neq\emptyset.$
If \FF\ satisfies the functional equation~(\ref{q:fe}), 
then \FF\ is nonzero if and only if $f_1(q)=1.$

For every positive integer $n$, let $\Omega(n)$ denote the number 
of not necessarily distinct prime factors of $n$.
If $n = p_1^{r_1}\cdots p_k^{r_k},$ then $\Omega(n) = r_1+\cdots + r_k.$

\bt
Let \pol\ be a nonzero sequence of polynomials that satisfies the functional
equation~(\ref{q:fe}).  The support of \FF\ is a prime semigroup.
If
\[
\supp(\FF) = S(P),
\]
where $P$ is a set of prime numbers, then the sequence \FF\ is 
completely determined by the set of polynomials
$\FF_P = \{f_p(q)\}_{p\in P}.$
\et

\pf 
Since \FF\ is nonzero, we have $f_1(q) = 1$ and so $1 \in \supp(\FF).$
If $m\in \supp(\FF)$ and $n \in \supp(\FF)$, then $f_m(q) \neq 0$ and $f_n(q)\neq 0,$
hence
\[
f_{mn}(q) = f_m(q)f_n(q^m)\neq 0
\]
and $mn\in \supp(\FF).$  Therefore, $\supp(\FF)$ is a semigroup.

Let $P$ be the set of prime numbers contained in $\supp(\FF)$.  Then $S(P) \subseteq \supp(\FF).$
If $n \in \supp(\FF)$ and the prime number $p$ divides $n$, then $n = pm$ for some 
positive integer $m$.  Since \FF\ satisfies the functional equation~(\ref{q:fe}),
we have
\[
f_n(q) = f_{pm}(q) = f_p(q)f_m(q^p) \neq 0,
\]
and so $f_p(q) \neq 0$, hence $p \in \supp(\FF)$ and $p \in P.$  
Since every prime divisor of $n$ belongs to $\supp(\FF)$,
it follows that $n \in S(P),$ and so $\supp(\FF) \subseteq S(P).$
Therefore, $\supp(\FF) = S(P)$ is a prime semigroup.

We use induction on $\Omega(n)$ for $n \in \supp(\FF)$
to show that the sequence $\FF_P = \{f_p(q)\}_{p\in P}$ determines \FF.
If $\Omega(n)=1,$ then $n = p \in P$ and $f_p(q) \in \FF_P.$
Suppose that $\FF_P$ determines $f_m(q)$ for 
all $m \in \supp(\FF)$ with $\Omega(m) \leq k$.
If $n \in \supp(\FF)$ and $\Omega(n) = k+1,$ then $n = pm,$ 
where $p \in P$, $m \in \supp(\FF)$, and $\Omega(m) = k.$ 
It follows that the polynomial $f_n(q) = f_p(q)f_m(q^p)$ 
is determined by $\FF_P.$
\eop

Let $P$ be a set of prime numbers, and let $S(P)$ be the semigroup generated by $P$.
Define the sequence \pol\ by
\[
f_n(q) = \left\{\ba{ll} 
[n]_q & \mbox{if $n \in S(P),$}\\
0 & \mbox{if $n \not\in S(P)$.}
\ea\right.
\]
Then \FF\ satisfies~(\ref{q:fe}) and $\supp(\FF) = S(P).$
Thus, every semigroup of the form $S(P)$ is the support of some 
sequence of polynomials satisfying the functional equation~(\ref{q:fe}).  

The following theorem provides a general method to construct
solutions of the functional equation~(\ref{q:fe}) with support $S(P)$
for any set $P$ of prime numbers.

\bt           \label{q:theorem:manyp}
Let $P$ be a set of prime numbers.  
For each $p \in P,$
let $h_p(q)$ be a nonzero polynomial such that 
\beq           \label{q:manyp0}
h_{p_1}(q) h_{p_2}(q^{p_1}) = h_{p_2}(q) h_{p_1}(q^{p_2})
\qquad\mbox{for all $p_1, p_2 \in P$.}
\eeq
Then there exists a unique sequence \pol\ with $\supp(\FF) = S(P)$
such that \FF\ satisfies the functional equation~(\ref{q:fe})
and $f_p(q) = h_p(q)$ for all $p \in P.$
\et

The proof uses three lemmas.

\bl                      \label{q:lemma:pk}
Let $p$ be a prime number and $h_p(q)$ a nonzero polynomial.
There exists a unique sequence of polynomials 
$\{f_{p^k}(q)\}_{k=0}^{\infty}$
such that $f_p(q) = h_p(q)$ and
\beq  \label{q:pk}
f_{p^{k}}(q) = f_{p^{i}}(q) f_{p^{j}}( q^{p^i}) 
\eeq
for all nonnegative integers $i,j$ and $k$
such that $i+j = k$.
\el

\pf
We define $f_1(q) = 1,$ $f_p(q) = h_p(q),$ and, by induction on $k$,
\beq          \label{q:defpk}
f_{p^{k}}(q) = f_p(q)f_{p^{k-1}}(q^p)
\eeq
for $k \geq 2.$
The proof of~(\ref{q:pk}) is by induction on $k.$  
Identity~(\ref{q:pk}) holds for $k = 0,1,$
and 2, and also for $i=0$ and all $j.$
Assume that~(\ref{q:pk}) is true for some $k \geq 1$, and 
let $k+1 = i+j,$ where $i \geq 1.$  
From the construction of the sequence $\{f_{p^k}(q)\}_{k=0}^{\infty}$ 
and the induction hypothesis we have
\bq
f_{p^{k+1}}(q) 
& = & f_{p}(q)f_{p^k}(q^p)  \\
& = & f_{p}(q) f_{p^{(i-1)+j}}(q^p) \\
& = & f_{p}(q) f_{p^{i-1}}(q^p) f_{p^j}( (q^p)^{p^{i-1}}) \\
& = & f_{p^{i}}(q) f_{p^{j}}( q^{p^i}).
\eq
Conversely, if the sequence $\{f_{p^k}(q)\}_{k=0}^{\infty}$
satisfies~(\ref{q:pk}), then, setting $i=1,$ we obtain~(\ref{q:defpk}),
and so the sequence $\{f_{p^k}(q)\}_{k=0}^{\infty}$
is unique.
\eop

\bl                      \label{q:lemma:p1p2}
Let $P = \{ p_1,p_2\},$ where $p_1$ and $p_2$ are distinct prime numbers, 
and let $S(P)$ be the semigroup generated by $P.$
Let $h_{p_1}(q)$ and $h_{p_2}(q)$ be nonzero polynomials such that
\beq           \label{q:h1h2}
h_{p_1}(q)h_{p_2}(q^{p_1}) = h_{p_2}(q)h_{p_1}(q^{p_2}).
\eeq
There exists a unique sequence of polynomials
$\{f_{n}(q)\}_{n\in S(P)}$ 
such that $f_{p_1}(q) = h_{p_1}(q),$ $f_{p_2}(q) = h_{p_2}(q),$
and
\beq  \label{q:p1p2}
f_{mn}(q) = f_m(q)f_n(q^m)\qquad\mbox{for all $m,n \in S(P)$.}
\eeq
\el

\pf
Every integer $n \in S(P)$ can be written uniquely in the form 
$n = p_1^ip_2^j$ for some nonnegative integers $i$ and $j$.
We apply Lemma~\ref{q:lemma:pk} to construct the sets of polynomials 
$\{f_{p_1^i}(q)\}_{i=0}^{\infty}$ and $\{f_{p_2^j}(q)\}_{j=0}^{\infty}$.
If $n = p_1^ip_2^j$ for positive integers $i$ and $j$, 
then we define
\beq         \label{q:defp1p2}
f_n(q) = f_{p_1^i}(q)f_{p_2^j}(q^{p^i}).
\eeq
This determines the set $\{f_{n}(q)\}_{n\in S(P)}$ .

We shall show that 
\beq             \label{q:12fe}
f_{p_1^i}(q)f_{p_2^j}(q^{p_1^i}) = f_{p_2^j}(q)f_{p_1^i}(q^{p_2^j})
\eeq
for all nonnegative integers $i$ and $j$.
This is true if $i=0$ or $j=0,$ so we can assume that $i\geq 1$ and $j\geq 1$.

The proof is by induction on $k=i+j$.
If $k = 2,$ then $i=j=1$ and the result follows from~(\ref{q:h1h2}).
Let $k \geq 2,$ and assume that equation~(\ref{q:12fe}) is true 
for all positive integers $i$ and $j$ such that $i+j\leq k.$
Let $i+j+1=k+1.$
By Lemma~\ref{q:lemma:pk} and the induction assumption,
\bq
f_{p_1^i}(q) f_{p_2^{j+1}}(q^{p_1^i}) 
& = & f_{p_1^{i}}(q) f_{p_2^j}( q^{p_1^i}) f_{p_2}(q^{p_1^i p_2^j})  \\
& = & f_{p_2^j}(q) f_{p_1^i}( q^{p_2^j})f_{p_2}(q^{p_2^jp_1^i})  \\
& = & f_{p_2^j}(q) f_{p_2}( q^{p_2^j})f_{p_1^i}(q^{p_2^{j+1}})  \\
& = & f_{p_2^{j+1}}(q) f_{p_1^i}(q^{p_2^{j+1}}).
\eq
Similarly,
\[
f_{p_1^{i+1}}(q) f_{p_2^j}(q^{p_1^{i+1}}) 
= f_{p_2^j}(q) f_{p_1^{i+1}}(q^{p_2^j}).
\]
This proves~(\ref{q:12fe}).

Let $m,n \in S(P).$  There exist nonnegative integers $i,j,k$, and $\ell$ such that
\[
m= p_1^ip_2^j\qquad\mbox{and}\qquad n= p_1^kp_2^{\ell}.
\]
Then
\bq
f_m(q)f_n(q^m)
& = & f_{p_1^i}(q)f_{p_2^j}(q^{p_1^i}) 
f_{p_1^k}(q^{p_1^ip_2^j})f_{p_2^{\ell}}(q^{p_1^{i+k}p_2^j}) \\
& = & f_{p_1^i}(q)f_{p_1^k}(q^{p_1^i})f_{p_2^j}(q^{p_1^{i+k}})f_{p_2^{\ell}}(q^{p_1^{i+k}p_2^j}) \\
& = & f_{p_1^k}(q)f_{p_1^i}(q^{p_1^k})f_{p_2^{\ell}}(q^{p_1^{i+k}})f_{p_2^j}(q^{p_1^{i+k}p_2^{\ell}}) \\
& = & f_{p_1^k}(q)f_{p_2^{\ell}}(q^{p_1^k})f_{p_1^i}(q^{p_1^kp_2^{\ell}})f_{p_2^j}(q^{p_1^{i+k}p_2^{\ell}}) \\
& = & f_n(q)f_m(q^n).
\eq
Setting $m = p_1^i$ and $n = p_2^j$ in~(\ref{q:p1p2}) gives~(\ref{q:defp1p2}),
and so the sequence of polynomials $\{f_{n}(q)\}_{n\in S(P)}$ is unique.
\eop

\bl                      \label{q:lemma:pr}
Let $P = \{ p_1,\ldots,p_r\}$ be a set consisting of $r$ prime numbers,
and let $S(P)$ be the semigroup generated by $P.$
 Let $h_{p_1}(q),\ldots,h_{p_r}(q)$ be nonzero polynomials such that
\beq           \label{q:hr}
h_{p_i}(q)h_{p_j}(q^{p_i}) = h_{p_j}(q)h_{p_i}(q^{p_j})
\eeq
for $i,j = 1,\ldots, r.$
There exists a unique sequence of polynomials
$\{f_{n}(q)\}_{n\in S(P)}$ 
such that $f_{p_i}(q) = h_{p_i}(q)$ for $i = 1,\ldots,r$, and
\beq  \label{q:pr}
f_{mn}(q) = f_m(q)f_n(q^m)\qquad\mbox{for all $m,n \in S(P)$.}
\eeq
\el

\pf
The proof is by induction on $r$.  
The result holds for $r=1$ by Lemma~\ref{q:lemma:pk}
and for $r=2$ by Lemma~\ref{q:lemma:p1p2}.
Let $r \geq 3,$ and assume that the Lemma holds for 
every set of $r-1$ primes.
Let $P' = P\setminus\{p_r\} = \{p_1,\ldots,p_{r-1}\}$.
By the induction hypothesis, there exists a unique sequence of polynomials
$\{f_{n}(q)\}_{n\in S(P')}$ 
such that $f_{p_i}(q) = h_{p_i}(q)$ for $i = 1,\ldots,r-1$, and
\[
f_{m'n'}(q) = f_{m'}(q)f_{n'}(q^{m'})\qquad\mbox{for all $m',n' \in S(P')$.}
\]
Every $n \in S(P)\setminus S(P')$ can be written uniquely
in the form $n = n'p_r^{a_r}$, where $n' \in S(P')$ and $a_r$
is a positive integer.
We define $f_{p_r^{a_r}}(q)$ by Lemma~\ref{q:lemma:pk} and
\beq         \label{q:defpr}
f_{n'p_r^{a_r}}(q) = f_{n'}(q)f_{p_r^{a_r}}(q^{n'})
\eeq
We begin by proving that
\beq        \label{q:kfe}
f_{n'}(q) f_{p_r^{a_r}}(q^{n'}) = 
f_{p_r^{a_r}}(q) f_{n'}(q^{p_r^{a_r}})
\eeq
for all $n' \in S(P')$ and $a_r \in \N.$

By Lemma~\ref{q:lemma:p1p2}, equation~(\ref{q:kfe}) is true if $n' = p_s^{a_s}$ 
for some prime $p_s \in P'$.
Let $n' = n''p_s^{a_s},$ where $n'' \in S\left(P\setminus\{p_s,p_r\}\right).$
By the induction assumption, 
\[
f_{n''}(q) f_{p_r^{a_r}}(q^{n''}) = 
f_{p_r^{a_r}}(q) f_{n''}(q^{p_r^{a_r}}),
\]
and so 
\bq
f_{n'}(q) f_{p_r^{a_r}}(q^{n'}) 
& = & f_{n''}(q)f_{p_s^{a_s}}(q^{n''})f_{p_r^{a_r}}(q^{n''p_s^{a_s}}) \\
& = & f_{n''}(q)f_{p_r^{a_r}}(q^{n''})f_{p_s^{a_s}}(q^{n''p_r^{a_r}}) \\
& = & f_{p_r^{a_r}}(q) f_{n''}(q^{p_r^{a_r}}) f_{p_s^{a_s}}(q^{n''p_r^{a_r}}) \\
& = & f_{p_r^{a_r}}(q) f_{n'}(q^{p_r^{a_r}}).
\eq
This proves~(\ref{q:kfe}).

Let $m,n \in S(P).$  We write $n = n'p_r^{a_r}$ and $m = m'p_r^{b_r}$,
where $m',n' \in S(P')$ and $a_r, b_r$ are nonnegative integers. 
Applying~(\ref{q:kfe}) and the induction assumption, we obtain
\bq
f_m(q)f_n(q^m)
& = & f_{m'}(q) f_{p_r^{b_r}}(q^{m'}) f_{n'}(q^{m'p_r^{b_r}}) f_{p_r^{a_r}}(q^{m'n'p_r^{b_r}}) \\
& = & f_{m'}(q) f_{n'}(q^{m'}) f_{p_r^{b_r}}(q^{m'n'}) f_{p_r^{a_r}}(q^{m'n'p_r^{b_r}}) \\
& = & f_{n'}(q) f_{m'}(q^{n'}) f_{p_r^{a_r}}(q^{m'n'}) f_{p_r^{b_r}}(q^{m'n'p_r^{a_r}}) \\
& = & f_{n'}(q) f_{p_r^{a_r}}(q^{n'}) f_{m'}(q^{n'p_r^{a_r}}) f_{p_r^{b_r}}(q^{m'n'p_r^{a_r}}) \\
& = & f_{n}(q) f_{m}(q^{n}).
\eq
This proves~(\ref{q:pr}).

Applying~(\ref{q:pr}) with $m = n'$ and $n = p_r^{a_r},$ we obtain~(\ref{q:defpr}).
This shows that the sequence $\{f_n(q)\}_{n\in S(P)}$ is unique,
and completes the proof of the Lemma.
\eop

We can now prove Theorem~\ref{q:theorem:manyp}.

\pf
If $P$ is a finite set of prime numbers, then we construct the set of
polynomials $\{f_n(q)\}_{n\in S(P)}$ by Lemma~\ref{q:lemma:pr}, and we define
$f_n(q) = 0$ for $n \not\in S(P).$  This determines the sequence \pol\ uniquely.

If $P$ is infinite, we write $P = \{p_i\}_{i=1}^{\infty}.$
For every positive integer $r$, let $P_r = \{p_i\}_{i=1}^{r}$
and apply Lemma~\ref{q:lemma:pr} to construct the set of polynomials 
$\{f_n(q)\}_{n\in S(P_r)}$.  
Since
\[
P_1 \subseteq \cdots \subseteq P_r \subseteq P_{r+1} \subseteq \cdots \subseteq P
\]
and 
\[
S(P_1) \subseteq \cdots \subseteq S(P_r) \subseteq S(P_{r+1}) \subseteq \cdots \subseteq S(P),
\]
we have
\[
\{f_n(q)\}_{n\in S(P_1)} \subseteq \cdots \subseteq
\{f_n(q)\}_{n\in S(P_r)} \subseteq \{f_n(q)\}_{n\in S(P_{r+1})} \subseteq \cdots.
\]
Define
\[
\{f_n(q)\}_{n\in S(P)} = \bigcup_{r=1}^{\infty} \{f_n(q)\}_{n\in S(P_{r})}
\]
Setting $f_n(q) = 0$ for all $n \not\in S(P)$ uniquely determines a sequence
\pol\ that satisfies\tfe\ and $f_p(q) = h_p(q)$ for all $p \in P.$
This completes the proof.
\eop

For example, for the set $P = \{2,5,7\}$, the reciprocal polynomials
\bq
h_2(q) & = & 1 - q + q^2 \\ 
h_5(q) & = & 1 - q + q^3 - q^4 + q^5 - q^7 + q^8 \\
h_7(q) & = & 1 - q + q^3 - q^4 + q^6 - q^8 + q^9 - q^{11} + q^{12}.
\eq 
satisfy the commutativity condition~(\ref{q:manyp0}). 
There is a unique sequence of polynomials \pol\ constructed 
from $\{h_2(q), h_5(q), h_7(q)\}$ by Theorem~\ref{q:theorem:manyp}.
Since
\[
f_p(q) = h_p(q) = \frac{[p]_{q^3}}{[p]_q} \qquad\text{for $p \in P = \{2,5,7\}$,}
\]
it follows that  
\[
f_n(q) = \frac{[n]_{q^3}}{[n]_q} \qquad\text{for all $n \in S(P)$.}
\]
We have $\deg(f_n) = 2(n-1)$ for all $n \in S(P).$

We can refine Problem~\ref{q:p1} as follows.

\bp        
Let $P$ be a set of prime numbers.  
Determine all polynomial sequences \pol\ with support $S(P)$ 
that satisfy the functional equation~(\ref{q:fe}).
\ep

\bp
Let $P$ and $P'$ be sets of prime numbers with $P \subseteq P'$,
and let \pol\ be a sequence of polynomials with support $S(P)$
that satisfies the functional equation~(\ref{q:fe}).
Under what conditions does there exist a sequence
$\mathcal{F}' = \{f'_n(q)\}_{n=1}^{\infty}$
with support $S(P')$ such that $\mathcal{F}'$ satisfies~(\ref{q:fe})
and $f'_p(q) = f_p(q)$ for all primes $p \in P$?
\ep

\bp
Let $S$ be a multiplicative subsemigroup of the positive integers. 
Determine all sequences $\{f_n(q)\}_{n\in S}$ of polynomials such that
\[
f_{mn}(q) = f_m(q)f_n(q^m) \qquad\mbox{for all $m,n \in S.$}
\]
\ep

This formulation of the problem of classifying solutions 
of the functional equation does not assume that $S$ is a semigroup
of the form $S = S(P)$ for some set $P$ of prime numbers.

\section{An arithmetic functional equation}
An {\em arithmetic function} is a function whose domain
is the set \N\ of natural numbers.
The {\em support} of the arithmetic function $\delta$ is
\[
\supp(\delta) = \{ n\in \N : \delta(n) \neq 0\}.
\] 

\bl                                \label{q:lemma:afe}
Let $S$ be a semigroup of the natural numbers,
and $\delta(n)$ a complex-valued arithmetic function 
that satisfies the functional equation
\beq   \label{q:afe}
\delta(mn) = \delta(m) + m\delta(n) \qquad\mbox{for all $m,n \in S.$}
\eeq
Then there exists a complex number $t$ such that
\[
\delta(n)  = t(n-1) \qquad\mbox{for all $n \in S$.}
\]
\el

\pf
Let $\delta(n)$ be a solution of the functional equation~(\ref{q:afe}) on $S$.
Setting $m=n=1$ in~(\ref{q:afe}), we obtain $\delta(1) = 0$.
For all $m,n \in S\setminus \{1\}$ we have
\[
\delta(m) + m\delta(n) = \delta(mn) = \delta(nm) = \delta(n) + n\delta(m),
\]
and so
\[
\frac{\delta(m)}{m-1} = \frac{\delta(n)}{n-1}.
\]
It follows that there exists a number $t$ such that
$\delta(n) = t(n-1)$ for all $n \in S$.
This completes the proof.
\eop

Note that if $\delta(n) = 0$ for some $n \in S\setminus \{1\},$ 
then $\delta(n) = 0$ for all $n \in S.$

Let $\deg(f)$ denote the degree of the polynomial $f(q)$.

\bl          \label{q:lemma:deg}
Let $\mathcal{F} = \{f_n(q)\}_{n=1}^{\infty}$
be a nonzero sequence of polynomials that satisfies
the functional equation~(\ref{q:fe}).
There exists a nonnegative rational number $t$ such that
\beq             \label{q:deg}
\deg(f_n) = t(n-1)\qquad \mbox{for all $n \in \supp(\FF)$.}
\eeq
\el

\pf
Let $S = \supp(\FF).$
The functional equation~(\ref{q:fe}) implies that
\[
\deg(f_{mn}) = \deg(f_m)+m\deg(f_n)\qquad\mbox{for all $m,n \in S,$}
\]
and so $\deg(f_n)$ is an arithmetic function on the semigroup $S$
that satisfies the arithmetic functional equation~(\ref{q:afe}).
Statement~(\ref{q:deg}) follows immediately from Lemma~\ref{q:lemma:afe}.
\eop

We note that, in Lemma~\ref{q:lemma:deg},
the number $t$ is rational but not necessarily integral.
For example, if $\supp(\FF) = \{ 7^k : k\in \N_0\}$ 
and 
\[
f_{7^k}(q) = q^{2(1+7+7^2+\cdots+7^{k-1})} = q^{(7^k-1)/3},
\]
then $t_1 = 1/3.$

An arithmetic function $\lambda(n)$ is
{\em completely multiplicative} if $\lambda(mn) = \lambda(m)\lambda(n)$
for all $m,n \in \N.$
A function $\lambda(n)$ is
{\em completely multiplicative on a semigroup $S$} 
if $\lambda(n)$ is a function defined on $S$
and $\lambda(mn) = \lambda(m)\lambda(n)$ for all $m,n \in S.$

\bt                \label{q:theorem:form}
Let $\mathcal{F} = \{f_n(q)\}_{n=1}^{\infty}$
be a nonzero sequence of polynomials that satisfies
the functional equation
\[
f_{mn}(q) = f_m(q)f_n(q^m).
\]
Then there exist a completely multiplicative arithmetic function
$\lambda(n)$, a nonnegative rational number $t$,
and a nonzero sequence \polg\ of polynomials such that 
\[
f_n(q) = \lambda(n)q^{t(n-1)}g_n(q)\qquad\mbox{for all $n \in \N$,}
\]
where
\benum
\item[(i)]
the sequence \GG\ satisfies the functional equation\fe,
\item[(ii)]
\[
\supp(\FF)=\supp(\GG)=\supp(\lambda),
\]
\item[(iii)]
\[
g_n(0)=1\qquad\mbox{for all $n \in \supp(\GG).$}
\]
\eenum
The number $t,$ the arithmetic function $\lambda(n)$, and the 
sequence \GG\ are unique.
\et

\pf
For every $n \in \supp(\FF)$ there exist a unique nonnegative integer $\delta(n)$
and polynomial $g'_n(q)$ such that $g_n'(0) \neq 0$ and
\[
f_n(q) = q^{\delta(n)}g'_n(q).
\]
Let $\lambda(n) = g'_n(0)$ be the constant term of $g'_n(q).$
Dividing $g_n'(q)$ by $\lambda(n),$ we can write
\[
g'_n(q) = \lambda(n)g_n(q),
\]
where $g_n(q)$ is a polynomial with constant term $g_n(0) = 1.$
Define $g_n(q) = 0$ and $\lambda(n) = 0$ 
for every positive integer $n \not\in \supp(\FF),$
and let \polg.  Then $\supp(\FF) = \supp(\GG) = \supp(\lambda).$
Since the sequence 
\[
\{\lambda(n)q^{\delta(n)}g_n(q)\}_{n=1}^{\infty}
\]
satisfies the functional equation, we have,
for all $m,n \in \supp(\FF)$,
\bq
\lambda(mn)q^{\delta(mn)}g_{mn}(q) 
& = & \lambda(m)q^{\delta(m)}g_m(q) \lambda(n)q^{m\delta(n)}g_n(q^m) \\
& = & \lambda(m)\lambda(n)q^{\delta(m)+m\delta(n)}g_m(q)g_n(q^m).
\eq
The polynomials $g_m(q), g_n(q),$ and $g_{mn}(q)$ have constant term 1, 
hence for all $m,n \in \supp(\FF)$ we have
\[
q^{\delta(mn)}= q^{\delta(m)+m\delta(n)},
\]
\[
\lambda(mn) = \lambda(m)\lambda(n),
\]
and
\[
g_{mn}(q) = g_m(q)g_n(q^m).
\]
It follows that $\lambda(n)$ is a completely multiplicative 
arithmetic function with $\supp(\FF),$
and the sequence $\{g_n(q)\}_{n=1}^{\infty}$ 
also satisfies the functional equation~({\ref{q:fe}).
Moreover,
\[
\delta(mn)= \delta(m)+m\delta(n) \qquad\mbox{for all $m,n \in\supp(\FF)$.}
\]
By Lemma~\ref{q:lemma:afe}, there exists a nonnegative rational number 
$t$ such that $\delta(n) = t(n-1)$.
This completes the proof.
\eop

\section{Classification problems}

Theorem~\ref{q:theorem:form} reduces the classification 
of solutions of the functional equation~(\ref{q:fe})
to the classification of sequences of polynomials \pol\ 
with constant term $f_n(0)= 1$ for all $n \in \supp(\FF).$

\bt          \label{q:theorem:new}
Let $\mathcal{F} = \{f_n(q)\}_{n=1}^{\infty}$
be a nonzero sequence of polynomials that
satisfies the functional equation~(\ref{q:fe}).
\benum
\item[(i)]
Let $\psi(q)$ be a polynomial such that $\psi(q)^m = \psi(q^m)$
for every integer $m \in \supp(\FF).$
Then the sequence $\{f_n(\psi(q))\}_{n=1}^{\infty}$
satisfies~(\ref{q:fe}). 
\item[(ii)] 
For every positive integer $t$,
the sequence $\{f_n(q^t)\}_{n=1}^{\infty}$
satisfies~(\ref{q:fe}). 
\item[(iii)]
The sequence of reciprocal polynomials 
$\{q^{\deg(f_n)}f_n(q^{-1})\}_{n=1}^{\infty}$ satisfies~(\ref{q:fe}).  
\eenum
\et

\pf
Suppose that $\psi(q)^m = \psi(q^m)$
for every integer $m \in \supp(\FF).$
Replacing $q$ by $\psi(q)$ in the polynomial identity\fe, 
we obtain
\[
f_{mn}(\psi(q)) = f_m(\psi(q))f_n(\psi(q)^m) = f_m(\psi(q))f_n(\psi(q^m))
\]
for all $m,n \in \supp(\FF).$  This proves~(i).

Since $(q^t)^m = (q^m)^t$ for all integers $t$, 
we obtain~(ii) from~(i) by choosing $\psi(q) = q^t.$

The reciprocal polynomial of $f(q)$ is
\[
\tilde{f}(q) = q^{\deg(f)}f(q^{-1}).
\]
Then
\bq
\tilde{f}_{mn}(q)
& = & q^{\deg(f_{mn})}f_{mn}(q^{-1}) \\
& = & q^{\deg(f_m)+m\deg(f_n)}f_{m}(q^{-1})f_n(q^{-m}) \\
& = & q^{\deg(f_m)} f_{m}(q^{-1}) q^{m\deg(f_n)} f_n\left(\left(q^{m}\right)^{-1}\right) \\
& = & \tilde{f}_m(q)\tilde{f}_n(q^m).
\eq
This proves~(iii).
\eop

For example, setting
\[
[n]_{q^t} = 1 + q^t + q^{2t} + \cdots + q^{(n-1)t},
\]
we see that $\{ [n]_{q^t}\}_{n=1}^{\infty}$ is a solution of~(\ref{q:fe})
with support \N.

The quantum integer $[n]_q$ is a self-reciprocal polynomial of $q$,
and $[n]_{q^t}$ is self-reciprocal for all positive integers $t$.
The reciprocal polynomial of the polynomial $q^{n-1}$ is 1.

The polynomials $\psi(q) = q^t$ are not the only polynomials 
that generate solutions of\tfe.  
For example, let $p$ be a prime number, 
and consider polynomials with coefficients
in the finite field $\Z/p\Z$ and solutions \pol\ of the functional equation
with $\supp(\FF) = S(\{p\}) = \{ p^k :k \in \N_0\}.$
Applying the Frobenius automorphism $z\mapsto z^p$,
we see that $\psi(q)^{m} = \psi(q^{m})$ for every polynomial $\psi(q)$
and every $m \in \supp(\FF).$

Here is another example of solutions of\fe\ generated by polynomials satisfying
$\psi(q)^{m} = \psi(q^{m})$ for $m \in \supp(\FF).$

\bt          \label{q:theorem:zeta}
Let $P$ be a nonempty set of prime numbers, 
and $S(P)$ the multiplicative semigroup generated by $P$.
Let $d$ be the greatest common divisor of the set $\{p-1: p \in P\}.$
For $\zeta \neq 0,$ let
\[
f_n(q) = \sum_{i=0}^{n-1} \zeta^iq^i = [n]_{\zeta q}\qquad\mbox{for $n \in S(P),$}
\]
and let $f_n(q)=0$ for $n\not\in S(P).$
If $\zeta$ is a $d$th root of unity,
then the sequence of polynomials 
\[
\mathcal{F} = \{f_n(q)\}_{n=1}^{\infty}
\]
satisfies the functional equation~(\ref{q:fe}).
Conversely, if \FF\ satisfies\fe, then $\zeta$ is a $d$th root of unity.
\et

\pf
Let $\zeta$ be a $d$th root of unity, and $\psi(q) = \zeta q.$
Since $p\equiv 1\pmod{d}$ for all $p\in P,$
it follows that $m\equiv 1\pmod{d}$ for all $m\in S(P).$
Therefore, if $m\in S(P),$ then
\[
\psi(q)^m = (\zeta q)^m = \zeta^mq^m = \zeta q^m = \psi(q^m).
\]
It follows from Theorem~\ref{q:theorem:new}
that the sequence of polynomials \pol, where
\[
f_n(q) = [n]_{\zeta q} = \sum_{i=0}^{n-1} \zeta^iq^i \qquad\mbox{for $n \in S(P)$}
\]
and $f_n(q) = 0$ for $n \not\in S(P)$, satisfies\tfe.

Conversely, suppose that \FF\ satisfies~(\ref{q:fe}).
Let $m,n \in S(P)\setminus\{ 1\}$.  
Since
\[
f_m(q)f_n(q^m) =  \left(\sum_{i=0}^{m-1} \zeta^iq^i\right)
\left( \sum_{j=0}^{n-1} \zeta^jq^{mj}\right) 
=  \sum_{i=0}^{m-1} \sum_{j=0}^{n-1} \zeta^{i+j}q^{i+mj},
\]
\[
f_{mn}(q) = \sum_{k=0}^{mn-1} \zeta^kq^k 
= \sum_{i=0}^{m-1} \sum_{j=0}^{n-1} \zeta^{i+mj}q^{i+mj},
\]
and
\[
f_{mn}(q) = f_m(q)f_n(q^m),
\]
it follows that
\[
\zeta^{i+j} = \zeta^{i+mj}
\]
for $0 \leq i \leq m-1$ and $0 \leq j \leq n-1.$ 
Then
\[
\zeta^{j(m-1)} = 1
\]
and
\[
\zeta^{m-1}=1 \qquad\mbox{for all $m \in S(P).$}
\]
Thus, $\zeta$ is a primitive $\ell$th root of unity for some positive
integer $\ell,$ and $\ell$ divides $m-1$ for all $m \in S(P).$
Therefore, $\ell$ divides $d,$ the greatest common divisor of the integers $m-1,$
and so $\zeta$ is a $d$th root of unity.
This completes the proof.
\eop

Let \pol\ and \polg\ be sequences of polynomials.
Define the product sequence
\[
\FF\GG = \{f_ng_n(q)\}_{n=1}^{\infty} 
\]
by $f_ng_n(q) = f_n(q)g_n(q).$

\bt
Let \FF\ and \GG\ be nonzero sequences of polynomials 
that satisfy the functional equation~(\ref{q:fe}).
The product sequence \FF\GG\ also satisfies~(\ref{q:fe}).  
Conversely, if $\supp(\FF) = \supp(\GG)$ and if \FF\ and
\FF\GG\ satisfy~(\ref{q:fe}), then \GG\ also satisfies~(\ref{q:fe}).  
The set of all solutions of~\tfe\ is an abelian semigroup,
and, for every prime semigroup $S(P)$,  the set $\Gamma(P)$ 
of all solutions \pol\ of\fe\ 
with $\supp(\FF) = S(P)$ is an abelian cancellation semigroup. 
\et

\pf
If \FF\ and \GG\ both satisfy~(\ref{q:fe}), then 
\bq
f_{mn}g_{mn}(q) & = & f_{mn}(q)g_{mn}(q) \\
& = & f_m(q)f_n(q^m)g_m(q)g_n(q^m) \\
& = & f_mg_m(q) f_ng_n(q^m),
\eq
and so \FF\GG\ satisfies\fe.
Conversely, if $m,n \in \supp(\FF)=\supp(\GG),$
\[
f_{mn}(q)g_{mn}(q) = f_m(q)g_m(q)f_n(q^m)g_n(q^m),
\]
and
\[
f_{mn}(q) = f_m(q)f_n(q^m),
\]
then
\[
g_{mn}(q) = g_m(q)g_n(q^m).
\]
Multiplication of sequences that satisfy\fe\ is associative and commutative.
For every prime semigroup $S(P)$, we define the sequence 
$\mathcal{I}_P =\{I_n(q)\}_{n=1}^{\infty}$
by $I_n(q) = 1$ for $n \in S(P)$ and $I_n(q) = 0$ for $n \not\in S(P)$.
Then $\mathcal{I}_P \in \Gamma(P)$ and 
$\mathcal{I}_P \FF = \FF$ for every $\FF \in \Gamma(P)$.
If $\FF, \GG, \HH \in \Gamma(P)$ and $\FF\GG = \FF\HH,$ 
then $\GG = \HH.$  Thus, $\Gamma(P)$ is a cancellation semigroup.
This completes the proof.
\eop

Let $S(P)$ be a prime semigroup, 
and let \pol\ and \polg\ be sequences of polynomials with support $S(P)$.
We define the sequence of rational functions $\FF/\GG$ by 
\[
\frac{\FF}{\GG} = \left\{ \frac{f_n}{g_n}(q)\right\}_{n=1}^{\infty},
\]
where
\[
\frac{f_n}{g_n}(q) = \frac{f_n(q)}{g_n(q)} \qquad \mbox{if $n \in S(P)$,} 
\]
and
\[
\frac{f_n}{g_n}(q) = 0 \qquad \mbox{if $n \not\in S(P)$.} 
\]
Then $\FF/\GG$ has support $S(P)$.
If \FF\ and \GG\ satisfy\tfe, then the sequence $\FF/\GG$ of rational functions
also satisfies\fe.

We recall the definition of the Grothendieck group of a semigroup.
If $\Gamma$ is an abelian cancellation semigroup, then there exists
an abelian group $K(\Gamma)$ and an injective semigroup homomorphism
$j: \Gamma \rightarrow K(\Gamma)$ such that if $G$ is any abelian group
and $\alpha$ a semigroup homomorphism from $\Gamma$ into $G$, then there 
there exists a unique group homomorphism $\tilde{\alpha}$ from $K(\Gamma)$
into $G$ such that $\tilde{\alpha} j = \alpha.$
The group $K(\Gamma)$ is called the {\em Grothendieck group}
of the semigroup $\Gamma.$

\bt
Let $S(P)$ be a prime semigroup, and let $\Gamma(P)$ be 
the cancellation semigroup of polynomial solutions of\tfe\ 
with support $S(P)$.
The Grothendieck group of $\Gamma(P)$ is the group of all sequences of rational functions
$\FF/\GG,$ where \FF\ and \GG\ are in $\Gamma(P).$
\et

\pf
The set $K(\Gamma(P))$ of all sequences of rational functions
of the form $\FF/\GG,$ where \FF\ and \GG\ are in $\Gamma(P),$
is an abelian group, and $\FF \mapsto \FF/\mathcal{I}_P$
is an imbedding of $\Gamma(P)$ into $K(\Gamma(P))$.
Let $\alpha:\Gamma(P) \rightarrow G$ be a homomorphism from $\Gamma(P)$
into a group $G$.  We define $\tilde{\alpha}:K(\Gamma(P)) \rightarrow G$
by
\[
\tilde{\alpha}\left(\frac{\FF}{\GG}\right) = 
\frac{\alpha(\FF)}{\alpha(\GG)}.
\]
If $\FF/\GG = \FF_1/\GG_1,$ then $\FF\GG_1 = \FF_1\GG.$
Since $\alpha$ is a semigroup homomorphism, we have
$\alpha(\FF)\alpha(\GG_1) = \alpha(\FF_1)\alpha(\GG),$
and so
\[
\frac{\alpha(\FF)}{\alpha(\GG)} = \frac{\alpha(\FF_1)}{\alpha(\GG_1)}.
\]
This proves that $\tilde{\alpha}:K(\Gamma(P)) \rightarrow G$
is a well-defined group homomorphism, and $\tilde{\alpha}j = \alpha.$
\eop

\bp
Does every sequence of rational functions
that satisfies\tfe\ and has support $S(P)$ 
belong to the group $K(\Gamma(P))$?
\ep

We recall that if \FF\ is a sequence of nonconstant polynomials
that satisfies~(\ref{q:fe}), then there exists a positive rational
number $t$ such that $\deg(f_n) = t(n-1)$
is a positive integer for all $n\in \supp(\FF).$
In particular, if $\supp(\FF) = \N$, or if $2 \in \supp(\FF)$, 
or, more generally, 
if $\{n-1:n\in \supp(\FF)\}$ is a set of relatively prime integers,
then $t$ is a positive integer.

The result below shows that the quantum integers are the unique solution 
of the functional equation\fe\ in the following important case.

\bt           \label{q:theorem:main}
Let \pol\ be a sequence of polynomials that satisfies 
the functional equation
\[
f_{mn}(q) = f_m(q)f_n(q^m)
\]
for all positive integers $m$ and $n$.
If $\deg(f_n) = n-1$ and $f_n(0) = 1$ 
for all positive integers $n$, then $f_n(q) = [n]_q$
for all $n$.
\et

Theorem~\ref{q:theorem:main} is a consequence of the following more general result.

\bt          \label{q:theorem:good}
Let \pol\ be a sequence of polynomials that
satisfies the functional equation
\[
f_{mn}(q) = f_m(q)f_n(q^m)
\]
for all positive integers $m$ and $n$.
If $\deg(f_n) = n-1$ and $f_n(0) = 1$ for all $n \in \supp(\FF),$
and if $\supp(\FF)$ contains 2 and some odd integer greater than 1,
then $f_n(q) = [n]_q$ for all $n\in \supp(\FF)$.
\et

\pf
Since $2 \in \supp(\FF)$, we have $\deg(f_2) = 1$ and $f_2(0) = 1$, 
hence
\[
f_2(q) = 1+aq
\]
for some $a \neq 0.$
If $n=2r+1 \geq 3$ is an odd integer in $\supp(\FF)$, then
\[
f_n(q) = 1 + \sum_{j=1}^{n-1}b_jq^j, \qquad\mbox{with $b_{n-1} \neq 0.$}
\]
We have
\bq
f_n(q)f_2(q^n) 
& = & \left(1+\sum_{j=1}^{n-1}b_jq^j\right)\left(1+aq^n\right) \\
& = & 1+\sum_{j=1}^{n-1}b_jq^j + aq^n + \sum_{j=1}^{n-1}ab_jq^{n+j}\\
& = & 1 + b_1q + b_2q^2 + \cdots + ab_1q^{n+1} + ab_2q^{n+2} + \cdots 
\eq
and 
\bq
f_2(q)f_n(q^2)
& = & (1+aq)\left(1 + \sum_{j=1}^{n-1}b_jq^{2j} \right) \\
& = & 1 + aq + \sum_{j=1}^{2r} b_j q^{2j} + \sum_{j=1}^{2r} ab_j q^{2j+1} \\
& = & 1 + aq + b_1q^2 + \cdots + b_{r+1}q^{n+1}+ ab_{r+1}q^{n+2} +\cdots.
\eq
The functional equation with $m=2$ gives
\beq               \label{q:2nfe}
f_n(q)f_2(q^n) = f_2(q)f_n(q^2).
\eeq
Equating coefficients in these polynomials, we obtain
\[
a = b_1 = b_2,
\]
\[
b_{r+1} = ab_1 = a^2,
\]
and
\[
ab_{r+1}=ab_2 = a^2.
\]
Since $a \neq 0,$ it follows that 
\[
a = 1
\]
and
\[
f_2(q)=1+q = [2]_q.
\]
By the functional equation, 
if $f_{2^{k-1}}(q) = [2^{k-1}]_q$ for some integer $k \geq 2,$
then
\bq
f_{2^k}(q) 
& = & f_{2^{k-1}}(q) f_2\left(q^{2^{k-1}}\right)  \\
& = & \left(1+q+q^2+\cdots + q^{2^{k-1}-1}\right) \left(1+q^{2^{k-1}}\right)  \\
& = & 1+q+q^2+\cdots + q^{2^{k}-1} \\
& = & [2^k]_q.
\eq
It follows by induction that $f_{2^k}(q) = [2^k]_q$ for all $k \in \N.$

Let $n = 2r+1$ be an odd integer in $\supp(\FF)$, $n \geq 3.$
Equation~(\ref{q:2nfe}) implies that
\[
1 + b_1q + \sum_{j=2}^{n-1} b_jq^j + q^n 
= 1 + q + \sum_{i=1}^{r} b_i \left( q^{2i} + q^{2i+1} \right),
\]
and so $1 = b_1 = b_r = b_{n-1}$ and
\[
b_i = b_{2i} = b_{2i+1}\qquad\mbox{for $i=1,\ldots,r-1$}.
\]
If $n = 3,$ then $b_1 = b_2 = 1$ and $f_3(q) = [3]_q.$
If $n = 5,$ then $r=2$ and $b_1 = b_2 = b_3 = b_4 = 1,$ hence $f_5(q) = [5]_q.$

For $n \geq 7$ we have $r \geq 3.$
If $1 \leq k \leq r-2$ and $b_i = 1$ for $i = 1,\ldots, 2k-1,$
then $k \leq 2k-1$ and so
\[ 
1 = b_k = b_{2k} = b_{2k+1}.
\]
It follows by induction on $k$ that $b_i = 1$ for $i = 1,\ldots, n-1,$
and $f_n(q) = [n]_q$ for every odd integer $n \in \supp(\FF)$.

If $2^kn \in \supp(\FF)$, where $n$ is odd, then
\[
f_{2^kn}(q) = f_{2^k}(q)f_n(q^{2^k}) = [2^k]_q[n]_{q^{2^k}} = [2^kn]_q.
\]
This completes the proof.
\eop

\bp
Let $t \geq 2$, and let \pol\ be a sequence of polynomials satisfying
the functional equation~(\ref{q:fe}) such that
$f_n(q)$ has degree $t(n-1)$ and $f_n(0)=1$ for all $n \in\N$.
Is \FF\ constructed from the quantum integers?
More precisely, do there exist positive integers 
$t_1,u_1,\ldots, t_k,u_k$ such that
\[
t = t_1u_1+\cdots+t_ku_k
\]
and, for all $n \in \N,$
\[
f_n(q) = \prod_{i=1}^k \left([n]_{q^{t_i}}\right)^{u_i}?
\]
\ep

\section{Addition of quantum integers}
It is natural to consider the analogous problem of addition of
quantum integers.  With the usual rule for addition of polynomials, 
$[m]_q + [n]_q \neq [m+n]_q$ for all positive integers
$m$ and $n$.  However, we observe that 
\[
[m]_q + q^m[n]_q = [m+n]_q 
\qquad\mbox{for all $m,n\in \N.$}
\]
This suggests the following definition.  Let \pol\ be a sequence of 
polynomials.  We define
\beq      \label{q:addfe}
f_m(q) \addq f_n(q) = f_m(q) + q^m f_n(q).
\eeq
If $h(q)$ is any polynomial, then the sequence \pol\ 
defined by $f_n(q) = h(q)[n]_q$
is a solution of the additive functional equation~(\ref{q:addfe}),
and, conversely, every solution of~(\ref{q:addfe}) is of this form.
This is discussed in Nathanson~\cite{nath02c}.


\begin{thebibliography}{99}
\bibitem{nath02c}  
M. B. Nathanson, 
Additive number theory and the ring of quantum integers, www.arXiv.org:
math.NT/0204006.                                                                                                                                                                                                                                               
\end{thebibliography}
\end{document}